\documentclass[12pt]{amsart}
\newcommand{\dsum}{\di\sum}
\renewcommand{\k}[1]{\ensuremath{\left({#1}\right)}}
\usepackage{multirow}\newcommand{\ds}{\dots}

\newcommand{\innn}{\ensuremath{\in\nn}}
\usepackage{centernot}
\usepackage{geometry,colortbl,tabularx}                
\newcommand{\p}{\ensuremath{\mathfrak{p}}}
\newcommand{\f}{\frac}

\newcommand{\z}{\zeta}
\newcommand{\cd}{\cdot}

\newcommand{\te}{\text}

\newcommand{\mb}{\mbox}

\newcommand{\h}{\texttt{h}}

\newcommand{\np}{\newpage}

\newcommand{\nn}{\mathbf{N}}

\usepackage{graphicx}
\usepackage{amssymb,enumerate}
\usepackage{amsmath,amsfonts,amsthm}
\usepackage{exscale,bm}
\usepackage{latexsym}
\usepackage{pict2e}\usepackage{verbatim}\usepackage[all,poly,curve,knot,arrow]{xy}
\usepackage{epstopdf}
\theoremstyle{plain}
\newtheorem{thm}{Theorem}[section]

\newtheorem{lem}[thm]{Lemma}

\newtheorem{prop}[thm]{Proposition}
\newtheorem{cor}[thm]{Corollary}

\theoremstyle{definition}
\newtheorem{ex}[thm]{Example}
\newtheorem{ques}[thm]{Question}

\newtheorem{defn}[thm]{Definition}
\newtheorem{rmk}[thm]{Remark}
\newtheorem{su}[thm]{Summary}
\newtheorem{ob}[thm]{Observation}

\newtheorem{fact}[thm]{Fact}

\newcommand{\sub}{\subseteq}

\newcommand{\be}{\begin{enumerate}}
\newcommand{\ee}{\end{enumerate}}
\newcommand{\bsu}{\begin{su}}
\newcommand{\esu}{\end{su}}

\newcommand{\bd}{\begin{defn}}
\newcommand{\ed}{\end{defn}}
\newcommand{\bp}{\begin{prop}}
\newcommand{\ep}{\end{prop}}
\newcommand{\bc}{\begin{cor}}
\newcommand{\ec}{\end{cor}}
\newcommand{\eq}{\end{ques}}
\newcommand{\bq}{\begin{ques}}
\newcommand{\bob}{\begin{ob}}
\newcommand{\eob}{\end{ob}}

\newcommand{\bl}{\begin{lem}}
\newcommand{\el}{\end{lem}}
\newcommand{\bt}{\begin{thm}}
\newcommand{\et}{\end{thm}}
\newcommand{\bpf}{\begin{proof}}
\newcommand{\epf}{\end{proof}}
\newcommand{\bex}{\begin{ex}}
\newcommand{\eex}{\end{ex}}

\newcommand{\di}{\displaystyle}
\newcommand{\cc}{\mathbb{C}}

\newcommand{\bft}{\begin{fact}}
\newcommand{\eft}{\end{fact}}
\newcommand{\brk}{\begin{rmk}}
\newcommand{\erk}{\end{rmk}}
\newcommand{\ba}{\begin{align*}}
\newcommand{\ea}{\end{align*}}
\newcommand{\tn}{\textnormal}

\newcommand{\bit}{\begin{itemize}}
\newcommand{\eit}{\end{itemize}}

\newcommand{\bcm}{}
\newcommand{\cref}[1]{(\ref{#1})}

\newcommand{\hf}{\hfill}

\newcommand{\eh}{\emph}

\newcommand{\lam}{\lambda}

\newcommand{\ff}[2]{\di\frac{#1}{#2}}

\newcommand{\led}[4]{$\xymatrix@1{{#1\,} \ar[r]^-{#2}_-{#3}& {\,#4}}$}
\newcommand{\wed}[4]{$\xymatrix@C=27pt@1{{#1\,} \ar@{->>}[r]^-{#2}_-{#3}& {\,#4}}$}

\begin{document}
\renewcommand{\h}{\hline}
\renewcommand{\arraystretch}{1.5}
\newcommand{\fib}{\text{fib\,}}
\newcommand{\fdp}{\text{FDP\,}}

\title{M\"obius functions of higher rank and Dirichlet series}
\author{Masato Kobayashi}
\date{\today}                                           
\thanks{masato210@gmail.com}

\address{Masato Kobayashi\\
Department of Engineering\\
Kanagawa University, 3-27-1 Rokkaku-bashi, Yokohama 221-8686, Japan.}


\maketitle
\begin{abstract} 
We introduce M\"obius functions of higher rank, a new class of arithmetic functions, so that the classical M\"obius function is of rank 2. With this idea, we evaluate Dirichlet series on the sum of reciprocal square of all $r$-free numbers. For the proof, Riemann zeta function and cyclotomic polynomials play a key role.
\end{abstract}
\tableofcontents
\np
\renewcommand{\lam}{\lambda}
\newcommand{\Set}[2]{\ensuremath{\left\{{#1}\,\middle|\,{#2}\right\}}}
\newcommand{\set}[1]{\ensuremath{\left\{{#1}\right\}}}
\renewcommand{\cd}{\cdots}

\newcommand{\mug}{\infty}
\newcommand{\q}{\quad}

\section{Introduction}

\subsection{Classical M\"obius and zeta functions}

The \eh{M\"obius function} plays an important role in number theory.
Its definition is simple:\,$\mu(n)=1$ and 
\[
\mu(n)=\begin{cases}
	(-1)^{k}&\te{$n=p_{1}\cd p_{k}$, primes $p_{j}$ all distinct,}\\
	0&\te{$p^{2}\,|\,n$ for some prime $p$.}\\
\end{cases}
\]

\eh{Riemann zeta function} is also another important topic in number theory. It is an analytic function of complex variable $s$  (pole at 1):
\[
\z(s)=\dsum_{n=1}^{\mug}\ff{1}{n^{s}}, \quad \text{Re}{(s)}>1.
\]
It has an infinite product (known as the \eh{Euler product}) 
expression:
\[
\z(s)=\di\prod_{p:\text{prime}}(1-p^{-s})^{-1}, 
\quad \text{Re}{(s)}>1.
\]
See Titchmarsh \cite{tit} for more details. Table \ref{zeven} shows the zeta values at positive even integers up to 20.

\newcommand{\lamn}{\lam(n)}
\renewcommand{\p}{\pi}
\newcommand{\hgt}[1]{\rule[0pt]{0pt}{#1}}
\newcommand{\dep}[1]{{\vrule width 0pt height 0pt depth #1}}
\newcommand{\Ome}{\Omega}
\renewcommand{\innn}{\in\nn}

\newcommand{\n}{\nu}
\newcommand{\mur}{\mu_{r}}

\newcommand{\dprod}{\di\prod}
\newcommand{\Om}{\Omega}

\newcommand{\murn}{\mur(n)}

{\renewcommand{\arraystretch}{2.3}
\begin{table}
\caption{zeta values at even positive integers}
\label{zeven}
\begin{center}
\begin{tabular}{c|c||c|ccccccccccccccccc}
$2n$&$\z(2n)$&$2n$&$\z(2n)$\\\h
2&$\ff{\p^{2}}{6}$&12&\dep{14pt}$\ff{691\p^{12}}{638512875}$\\\h
4&$\ff{\p^{4}}{90}$&14&\dep{14pt}$\ff{2\p^{14}}{18243225}$\\\h
6&$\ff{\p^{6}}{945}$&16&\dep{14pt}$\ff{3617\p^{16}}{325641566250}$\\\h
8&$\ff{\p^{8}}{9450}$&18&\dep{14pt}$\ff{43867\p^{18}}{38979295480125}$\\\h
10&$\ff{\p^{10}}{93555}$&20&\dep{14pt}$\ff{174611\p^{20}}{1531329465290625}$
\end{tabular}
\end{center}
\end{table}%
}

There is a deep relation between the M\"obius and zeta functions;
we can ``invert" $\zeta(s)$:
\[
\ff{1}{\zeta(s)}=\dsum_{n=1}^{\mug}\ff{\mu(n)}{n^{s}}, \q \textnormal{Re}{(s)}>1.
\]
For instance, when $s=2$, we obtain the inverse of Euler's work $\z(2)=\p^{2}/6$ as 
\[
\k{1+\ff{1}{2^{2}}+\ff{1}{3^{2}}+\ff{1}{4^{2}}+\ff{1}{5^{2}}+\cd}^{-1}
=
\ff{1}{\zeta(2)}=
\dsum_{n=1}^{\mug}\ff{\mu(n)}{n^{2}}
=1-\ff{1}{2^{2}}-\ff{1}{3^{2}}-\ff{1}{5^{2}}+\cd. 
\]

\subsection{Main results}

In this article, we introduce M\"obius functions of higher rank, a new class of arithmetic functions, so that the classical M\"obius function is of rank 2.  

\begin{center}
\begin{minipage}[c][100pt]{300pt}
\xymatrix@=9mm{
*+[F]{\text{\,\strut classical M\"{o}bius function $\mu=\mu_{2}$}\,}
\ar@{->}_-*\txt{}[d]\\
*+[F]{\text{\,\strut M\"{o}bius functions of higher rank $\mu_{r}$
 ($r=1, 2, \dots, \mug$)\,
}}\\
}
\end{minipage}
\end{center}

For a positive integer $r$, say $n$ is \eh{$r$-free} if 
there exists some prime $p$ such that $p^{r}$ divides $n$;
thus, 2-free is square-free and 3-free is cube-free as usually said. 
We will see that $\mu_{r}$ is similar to $\mu=\mu_{2}$: $\mu_{r}(n)\ne 0$ if and only if $n$ is $r$-free (Section \ref{sec3}).

The main result of this article is to evaluate several Dirichlet series
\[
\sum_{
n:r\textnormal{-free}}
\ff{\mu_{r}(n)}{n^{s}}
\]
with $r\in \{3, 4, 5\}$ and $s\in\{2, 3\}$. 

\begin{thm}[$s=2$]
The following equalities hold:
\begin{quote}
\begin{enumerate}
	\item 
		$\di\sum_{
	\substack
	{n:
	\textnormal{3-free}
	}
	}
	\ff{\mu_{3}(n)}{n^{2}}=\ff{45045}{691\pi^{4}}.
	$
\item 
		$\di\sum_{
	\substack
	{n:
	\textnormal{4-free}
	}
	}
	\ff{\mu_{4}(n)}{n^{2}}=\ff{630}{\pi^{6}}.
	$
		\item 
$	\di\sum_{
	\substack
	{n:
	\textnormal{5-free}
	}
	}
	\ff{\mu_{5}(n)}{n^{2}}=\ff{1091215125}{174611\p^{8}}
$.
\end{enumerate}
\end{quote}
\end{thm}
\begin{thm}[$s=3$]
\[
		\k{\di\sum_{
	\substack
	{n:
	\textnormal{3-free}
	}
	}
	\ff{1}{n^{3}}
}
	\k{\di\sum_{
	\substack
	{n:
	\textnormal{3-free}
	}
	}
	\ff{\mu_{3}(n)}{n^{3}}}
	=\ff{41247931725}{43867\p^{12}}.
	\]
\end{thm}

Here, $\mu_{3}, \mu_{4}, \mu_{5}$ are the M\"obius functions of rank 3, 4, 5 respectively. For the proofs, it is key to understand interactions of the following three concepts:
\begin{quote}
	\begin{itemize}
		\item Euler product for Riemann zeta function
		\item M\"obius functions of higher rank
		\item Cyclotomic polynomials
	\end{itemize}
\end{quote}
We will give these details later.

Additional results: it is possible to generalize the ``M\"obius inversion formula" to higher rank: for $n=p_{1}^{m_{1}}\cd p_{k}^{m_{k}}$, the prime factorization of 
$n$, let 
\[
m_{*}(n)=|\{j\,|\,m_{j}\equiv 3, \tn{4 mod 6}\}|.
\]

\begin{thm}[M\"obius inversion of rank 3]
\[
\k{\dsum_{n=1}^{\mug}\ff{\mu_{3}(n)}{n^{s}}}
^{-1}=
\sum_{
\substack
{n=p_{1}^{m_{1}}\cd p_{k}^{m_{k}}\\
m_{j}\not\equiv 2, \tn{5 mod }6
}
}
\ff{(-1)^{m_*(n)}}{n^{s}}.
\]
\end{thm}
As a by-product, we get a new expression of $\pi$:
\begin{align*}
	\p&=
	\k{\ff{45045}{691}
	\k{1+\ff{1}{2^{2}}+\ff{1}{3^{2}}+\ff{1}{5^{2}}+\ff{1}{6^{2}}
+\ff{1}{7^{2}}-\ff{1}{8^{2}}+\ff{1}{10^{2}}+\cd}}^{1/4}.
\end{align*}



\subsection{Notation}

	\begin{itemize}
	\item Let $\nn$ denote the set of positive integers.
	In addition, $\nn^{2}$ means 
the set of square numbers $\{1^{2}, 2^{2}, 3^{2}, \ds\}$.
	\item Often, writing 
	\[
n= p_{1}^{m_{1}}\cd p_{k}^{m_{k}}
\]
means the factorization of $n$ into \eh{distinct} prime numbers  ($p_{i}\ne p_{j}$ for $i\ne j$) with each $m_{j}$ positive unless otherwise specified.
		\item $d\,|\,n$ means $d$ divides $n$.
	\item $\dprod_{p}$ indicates an infinite product over all primes $p$.
	\end{itemize}

\section{Preliminaries}

Let us begin with 
recalling some fundamental definitions and facts on 
arithmetic functions; you can find this topic in a standard textbook on number theory as Apostol \cite{apostol}.
We thus omit most of the proofs here.

\subsection{Arithmetic functions}

An \eh{arithmetic function} is a map 
\[
f:\nn\to\cc.\]
\begin{ex}\hf
\begin{itemize}
	\item {\bf M\"obius function:} 
\[
\mu(n)=\begin{cases}
	1&n=1,\\
	(-1)^{k}&\te{$n=p_{1}\cd p_{k}$, primes $p_{j}$ all distinct,}\\
	0&\te{$p^{2}\,|\,n$ for some prime $p$.}\\
\end{cases}
\]
	\item{\bf Omega function:} 
	For $n=p_{1}^{m_{1}}\cd p_{k}^{m_{k}}$ with $p_{j}$ primes, $$\Ome(n)=m_{1}+\cd+m_{k}.$$
	\item {\bf Liouville function:} $\lam(n)=(-1)^{\Ome(n)}$.
	\item {\bf Characteristic function:} 
	For a subset $A\sub \nn$, 
	\[\chi_{A}(n)=\begin{cases}
	1&n\in A,\\
	0&n\not\in A.
\end{cases}
\]
In particular, $|\mu(n)|$ is a characteristic function of the set of 2-free numbers.
	\item {\bf Constant function:} $1(n)=1$ for all $n$.
	\item {\bf unit function:}
	$u(n)=\begin{cases}
	1&n=1,\\
	0&n\ne1.
	\end{cases}
	$
\end{itemize}
\end{ex}

We say that an arithmetic function $f$ is \eh{multiplicative} if $f(1)=1$ and 
\[
f(mn)=f(m)f(n)
\text{ \q whenever \,\, $\gcd(m, n)=1$.}
\]
It is easy to check that $\mu, \lam, 1, u$ are all multiplicative.

{\renewcommand{\arraystretch}{1.5}
\begin{table}
\caption{Arithmetic functions}
\begin{center}
\begin{tabular}{c|ccccccccccccccccccc}
$n$&1&2&3&4&5&6&7&8&9&10&$\cd$\\\hline
$\mu(n)$&1&$-1$&$-1$&0&$-1$&1&$-1$&0&0&1&$\cd$\\
$\Omega(n)$&0&1&1&2&1&2&1&3&2&2&$\cd$\\
$\lam(n)$&1&$-1$&$-1$&1&$-1$&1&$-1$&$-1$&1&1&$\cd$\\
$1(n)$&1&1&1&1&1&1&1&1&1&1&$\cd$\\
$u(n)$&1&0&0&0&0&0&0&0&0&0&$\cd$\\
\end{tabular}
\end{center}
\end{table}%
}

\subsection{Dirichlet series}

For two arithmetic functions $f$ and $g$, define 
the \eh{Dirichlet product} $f*g$ by
\[
(f*g)(n)=\dsum_{d\,|\,n}f(d)g\k{\ff{n}{d}}.
\]
The unit function $u$ satisfies 
\[
f*u=u*f=f
\]
for all arithmetic functions $f$.
If $f*g=g*f=u$, then we write $g=f^{-1}$ and call it the \eh{Dirichlet inverse} of $f$; assuming $f(1)\ne0$, there exists $f^{-1}$.

\begin{fact}
Let $f$ and $g$ be arithmetic functions.
Suppose they are multiplicative.
Then, so are $f*g$ and $f^{-1}$.
\end{fact}
In this way, multiplicative functions form a group and $u$ is indeed a  group-theoretic unit.

\begin{rmk}
If multiplicative functions $f, g$ satisfy
\begin{quote}
$f(p^{m})=g(p^{m})$ \,\, for all primes $p$ and $m\ge 1$,
\end{quote}
then $f(n)=g(n)$ for all $n\innn$.
Hence, to determine a multiplicative function, it is enough to know 
values only at prime powers.
\end{rmk}

A \eh{Dirichlet series} for $f$ is a series in the form
\[
\dsum_{n=1}^{\mug}\ff{f(n)}{n^{s}}
\]
for a complex number $s$ (in this article, we deal with only $s=2, 3$ and convergent series).
Riemann zeta function is an example of such series with $f(n)=1(n)=1$.

Observe that 
\[
\k{\sum_{n=1}^{\mug}\ff{f(n)}{n^{s}}}
\k{\sum_{n=1}^{\mug}\ff{g(n)}{n^{s}}}
=\sum_{n=1}^{\mug}\ff{(f*g)(n)}{n^{s}}
\]
for all $f, g$.
Then classical results 
\[
(\mu*1)(n)=
\dsum_{d\,|\,n}\mu(d)=
\begin{cases}
1&n=1,\\
0&n\ne1
\end{cases}
\]
and 
\[
(\lam*1)(n)=\chi_{\nn^{2}}(n)=
\begin{cases}
	1&n=N^{2} \text{ for some } N,\\
	0&\text{otherwise}
\end{cases}
\]
imply the following:
\begin{fact}\hfill
\begin{quote}
\begin{enumerate}
	\item 
$
\k{\dsum_{n=1}^{\mug}\ff{\mu(n)}{n^{s}}}
\k{\dsum_{n=1}^{\mug}\ff{1}{n^{s}}}
=1.
$\dep{20pt}
	\item 
$\k{\dsum_{n=1}^{\mug}\ff{\lam(n)}{n^{s}}}
\k{\di\dsum_{n=1}^{\mug}\ff{1}{n^{s}}}
=\dsum_{n:\textnormal{\,square}}\ff{1}{n^{s}}=
\sum_{N=1}^{\mug}\ff{1}{(N^{2})^{s}}=\z(2s).
$
\end{enumerate}
\end{quote}
\end{fact}

As a consequence, 
when $s=2$, we have
\[
\dsum_{n=1}^{\mug}\ff{\mu(n)}{n^{2}}=
\ff{1}{\z(2)}=\ff{6}{\p^{2}} 
\q \text{and} \q 
\dsum_{n=1}^{\mug}\ff{\lamn}{n^{2}}=
\ff{\z(4)}{\zeta(2)}=\ff{\p^{2}}{15}.
\]
Once we introduce the M\"obius functions of higher rank $\mu_{r}$ in the next section, we can regard these as extremal cases at $r=2$ and $r=\mug$ (as shown in Table \ref{mres}):
\[
\dsum_{n:\tn{2-free}}\ff{\mu_{2}(n)}{n^{2}}=\ff{6}{\p^{2}} 
\q \tn{and} \q 
\dsum_{n:\mug\tn{-free}}\ff{\mu_{\mug}(n)}{n^{2}}=
\ff{\p^{2}}{15}.
\]

{\renewcommand{\arraystretch}{2.5}
\begin{table}
\caption{Main results $(s=2)$}
\label{mres}
\begin{center}
	\begin{tabular}{c|ccccccc}\h
	series&	value&zeta expression&factor of Euler product\\\h
\dep{20pt}$\dsum_{
\substack
{n:\textnormal{2-free}
}
}\ff{\mu_{2}(n)}{n^{2}}$
	&	$\ff{6}{\p^{2}}$&$\ff{1}{\z(2)}$&$1-p^{-2}$\\\h
\dep{20pt}$\dsum_{
\substack
{n:\textnormal{3-free}
}
}\ff{\mu_{3}(n)}{n^{2}}$
	&	$\ff{45045}{691\p^{4}}$&$\ff{\z(4)\z(6)}{\z(2)\z(12)}$&$1-p^{-2}+p^{-4}$\\\h
\dep{20pt}$\dsum_{
\substack
{n:\textnormal{4-free}
}
}\ff{\mu_{4}(n)}{n^{2}}$
	&	$\ff{630}{\p^{6}}$&$\ff{\z(4)}{\z(2)\z(8)}$&$1-p^{-2}+p^{-4}-p^{-6}$\\\h
\dep{20pt}$\dsum_{
\substack
{n:\textnormal{5-free}
}
}\ff{\mu_{5}(n)}{n^{2}}$
	&	$\ff{1091215125}{174611\p^{8}}$&$\ff{\z(4)\z(10)}{\z(2)\z(20)}$&$1-p^{-2}+p^{-4}-p^{-6}+p^{-8}$\\\h
$\cd$&&$\cd$&	\\\h
\dep{20pt}	$\dsum_{
\substack
{n:\textnormal{$\mug$-free}
}
}\ff{\mu_{\mug}(n)}{n^{2}}$
	&	$\ff{\p^{2}}{15}$&$\ff{\z(4)}{\z(2)}$&$1-p^{-2}+p^{-4}-p^{-6}+p^{-8}-\cd$\\\h
\end{tabular}
\end{center}
\end{table}}



\section{M\"obius functions of higher rank}
\label{sec3}

For each natural number $r$ or $``r=\mug"$,
define an arithmetic function
\[
\mu_{r}:\nn\to\{-1, 0, 1\}\]
by
\[
\mu_{r}(n)=
\begin{cases}
	1&n=1,\\
	(-1)^{m_{1}+\cd+m_{k}}&
	n=p_{1}^{m_{1}}\cd p_{k}^{m_{k}}, \text{all } m_{j}< r,\\
	0&\text{$p^{r}\,|\,n$ for some prime $p$}.\\
\end{cases}
\]
For $r=\mug$, we understand that 
$m_{j}< \mug$ always holds and $p^{\mug}\,|\,n$ never happens.

\begin{ex}
$r=1$: This is just the unit function.
\[
\mu_{1}(n)=u(n)=\begin{cases}
	1&n=1,\\
	0&n\ne 1 \text{\,\,\,(that is, $p^{1}|n$ for some prime $p$)}.
\end{cases}
\]
$r=2$: the classical M\"obius function.
\[
\mu_{2}(n)=\mu(n)=\begin{cases}
	(-1)^{k}&\te{$n=p_{1}p_{2}\cd p_{k}$},\\
	0&\text{$p^{2}\,|\,n$ for some prime $p$}.
\end{cases}
\]
$r=3$:
\[
\mu_{3}(n)=\begin{cases}
	(-1)^{m_{1}+\cd+m_{k}}&\te{$n=p_{1}^{m_{1}}\cd p_{k}^{m_{k}}$, $m_{j}<3$},\\
	0&\text{$p^{3}\,|\,n$ for some prime $p$}.
\end{cases}
\]
$r=\mug$: the Liouville function.
\[
\mu_{\mug}(n)=\lam(n)=
(-1)^{m_{1}+\cd+m_{k}}\q \te{$n= p_{1}^{m_{1}}\cd p_{k}^{m_{k}}$}.
\]
\end{ex}

\begin{defn}
All together, we call $\{\mu_{r}\}_{r=1}^{\mug}$ the \eh{M\"obius functions of higher rank}.
\end{defn}


It follows by definition 
\[
\mu_{r}(p^{m})=\begin{cases}
	(-1)^{m}&m<r,\\
	0&m\ge r
\end{cases}
\]
for a prime $p$ and $m\ge1$.
Observe that each $\mu_{r}$ is multiplicative;
in particular, $|\mu_{r}|$ is a characteristic function of $r$-free numbers.

We have already seen that 
\[
\mu*1=u \text{\q and \q} \lam*1=\chi_{\nn^{2}}.
\]
Now understand this 
as $\mu_{2}*1=u$ and $\mu_{\mug}*1=\chi_{\nn^{2}}$.
A natural question is: what is $\mu_{r}*1$ for $3 \le r<\mug$?
Since $\mu_{r}$ and $1$ are both multiplicative, 
so is $\mu_{r}*1$. 
Now let us see what $(\mu_{r}*1)(p^{m})$ is.
\begin{prop}
Let $r\ge 3$ and $m\ge 1$.\\
If $m<r$, then 
\[
(\mur*1)(p^{m})=
\begin{cases}
	1& \textnormal{$m$ even},\\
	0& \textnormal{$m$ odd}.
\end{cases}
\]
\renewcommand{\tn}{\textnormal}
If $m\ge r$, then 
\[
(\mur*1)(p^{m})=\begin{cases}
	1&\tn{$r$ odd},\\
	0&\tn{$r$ even}.
\end{cases}
\]\end{prop}
\begin{proof}
Suppose $m<r$. Then 
\begin{align*}
	(\mur*1)(p^{m})&=\sum_{d\,|\,p^{m}}\mu_{r}(d)
	\\&=\mu_{r}(1)+\mur(p)+\mur(p^{2})+\cd+\mur(p^{m})
\\&=1+(-1)+1+\cd+(-1)^{m}
\\&=\begin{cases}
	1& \tn{$m$ even},\\
	0& \tn{$m$ odd}.
\end{cases}
\end{align*}
If $m\ge r$, then 
\begin{align*}
	(\mur*1)(p^{m})&=\sum_{d\,|\,p^{m}}\mu_{r}(d)
	\\&=\mu_{r}(1)+\mur(p)+\mur(p^{2})+\cd+\mur(p^{m})
	\\&=\mu_{r}(1)+\mur(p)+\mur(p^{2})+\cd+\mur(p^{r-1})+0+\cd+0
	\\&=1+(-1)+1+\cd+(-1)^{r-1}
	\\&=\begin{cases}
	1&\tn{$r$ odd},\\
	0&\tn{$r$ even}.
\end{cases}
\end{align*}
\end{proof}

Consequently, for $n= p_{1}^{m_{1}}\cd p_{k}^{m_{k}}$, the integer
\[
(\mu_{r}*1)(n)=(\mur*1)(p_{1}^{m_{1}}\cd p_{k}^{m_{k}})
=(\mur*1)(p_{1}^{m_{1}})\cd (\mur*1)(p_{k}^{m_{k}})
\]
is 1 if and only if all of 
factors $(\mur*1)(p_{j}^{m_{j}})$ are 1.
Otherwise, i.e., $(\mur*1)(p_{j}^{m_{j}})=0$ for some $j$, 
it is $0$.
This naturally leads to an interpretation of 
$\mu_{r}*1$ as a characteristic function of some set as follows. 
For each $r\ge3$, define $M_{r}$, a subset of $\nn$:
\begin{itemize}
	\item $r$ odd or $r=\mug$: square numbers.
	\[
M_{3}=M_{5}=\cd=M_{\mug}=
=\nn^{2}\,(=
\Set{n\innn}{n=p_{1}^{m_{1}}\cd p_{k}^{m_{k}}, m_{j}
\text{ all even}
}).
\]
	\item $r$ even: ranked square numbers.
	\[
M_{r}=\Set{n\innn}{n=p_{1}^{m_{1}}\cd p_{k}^{m_{k}}, m_{j}
\text{ all even}, m_{j}<r
}.
\]
\end{itemize}
The sets $M_{r}$'s ($r$ even) are increasing:
\[
M_{4}\subset M_{6}\subset M_{8}\subset\cd\subset M_{\mug}=\nn^{2}.
\]

\begin{ex}
\begin{align*}
	M_{4}&=\Set{n\in\nn}{n=p_{1}^{m_{1}}\cd p_{k}^{m_{k}}, m_{j}
\text{ all even}, m_{j}<4
}
	\\&=
	\Set{n\in\nn}{n=1 \mb{ or }n=p_{1}^{2}\cd p_{k}^{2}}
	\\&=\{1, 4, 9, 25, 36, 49, 100, 121, 169, \ds\}.
\end{align*}
\end{ex}


\begin{prop}
Let $M_{1}=\nn$ and $M_{2}=\{1\}$. Then, for each $r\in \nn\cup\{\mug\}$, the Dirichlet product $\mu_{r}*1$ is a characteristic function of the set $M_{r}$:
\[
(\mu_{r}*1)(n)=
\begin{cases}
	1&n\in M_{r},\\
	0&n\not\in M_{r}.
\end{cases}
\]
\end{prop}


{\renewcommand{\arraystretch}{1.25}
\begin{table}
\caption{$\mu_{3}$ and $\mu_{3}*1$}
\begin{center}
\begin{tabular}{c|ccccccccccccccccccc}
$n$&1&2&3&4&5&6&7&8&9&10&$\cd$\\\hline
$\mu_{3}(n)$&$1$&$-1$&$-1$&1&$-1$&1&$-1$&0&1&1&$\cd$\\
$(\mu_{3}*1)(n)$&1&0&0&1&0&0&0&0&1&0&$\cd$\\
\end{tabular}
\end{center}
\end{table}%
}

\begin{prop}[M\"obius functions of higher rank and zeta]
For $r\ge 1$ and $\te{Re}{(s)}>1$, we have 
\[
\k{\dsum_{n=1}^{\mug}\ff{\mu_{r}(n)}{n^{s}}}
\zeta(s)
=\sum_{n\in M_{r}}\ff{1}{n^{s}}.
\]
\end{prop}

\begin{proof}
This statement is equivalent to $\mu_{r}*1=\chi_{M_{r}}$.
\end{proof}

For clarity, we sometimes prefer to write 
\[
\dsum_{n=1}^{\mug}\ff{\mu_{r}(n)}{n^{s}}
=
\dsum_{n:\textnormal{$r$-free}}\ff{\mur(n)}{n^{s}}.
\]
In the next section, we will compute such sums for $s=2$.

\section{Main Theorems}

Before going into main theorems, we briefly recall an important family of polynomials in number theory for convenience.

\subsection{Cyclotomic polynomials}

The \eh{cyclotomic polynomial} for $n$ is 
\[
\Phi_{n}(x)=
\prod_{
\substack{1\le k\le n\\
\gcd(k, n)=1}
}(x-e^{2\pi i n/k}).
\]
This is indeed a polynomial of integer coefficients.
\begin{ex}
\[
\Phi_{1}(x)=x-1, \q
\Phi_{2}(x)=x+1, \q \textnormal{ and\q} 
\Phi_{3}(x)=x^{2}+x+1.
\]
\end{ex}
An important relation to the M\"obius function is:
\begin{fact}
\[
\Phi_{n}(x)=\prod_{d\,|\,n}(x^{d}-1)^{\mu\k{\f{n}{d}}.
}
\]
\end{fact}
Exponents are $0, \pm 1$
so that $\Phi_{n}(x)$ (and $\Phi_{n}(x)^{-1}$ also) is a product of $(x^{d}-1)$'s. Note that a factor $x^{d}-1$ looks like ``$1-p^{-s}$" in the 
Euler product of $\z(s)$; this idea will play a key role in the proofs below.

\subsection{Theorem ($s=2$)}

We are now ready for computing three series in the middle of Table \ref{mres}.

\begin{thm}\label{mth1}
\hf
\begin{quote}
\begin{enumerate}
	\item 
		$\di\sum_{
	\substack
	{n:
	\textnormal{3-free}
	}
	}
	\ff{\mu_{3}(n)}{n^{2}}=\ff{45045}{691\pi^{4}}.
	$
\item 
		$\di\sum_{
	\substack
	{n:
	\textnormal{4-free}
	}
	}
	\ff{\mu_{4}(n)}{n^{2}}=\ff{630}{\pi^{6}}.
	$
		\item 
$	\di\sum_{
	\substack
	{n:
	\textnormal{5-free}
	}
	}
	\ff{\mu_{5}(n)}{n^{2}}=\ff{1091215125}{174611\p^{8}}
$.
\end{enumerate}
\end{quote}
\end{thm}


\newcommand{\vi}{\\[.1in]}\newcommand{\vii}{\\[.2in]}

\begin{proof}[Proof of \tn{(1)}]
Note that
\begin{align*}
	\Phi_{12}(x)&=
	(x-1)^{\mu(12)}
	(x^{2}-1)^{\mu(6)}
	(x^{3}-1)^{\mu(4)}
	(x^{4}-1)^{\mu(3)}
	(x^{6}-1)^{\mu(2)}
	(x^{12}-1)^{\mu(1)}
	\\&=\ff{(1-x^{2})(1-x^{12})}{(1-x^{4})(1-x^{6})}
	=1-x^{2}+x^{4}.
\end{align*}
Then, we have 
\begin{align*}
	\dsum_{n:\text{3-free}}\ff{\mu_{3}(n)}{n^{2}}&=
	1+\sum_{
\substack
{0<m_{j}<3
}
}
\ff{(-1)^{m_{1}+\cd+m_{k}}}{(p_{1}^{m_{1}}\cd p_{k}^{m_{k}})^{2}}
	\\&=\dprod_{p}
(1-p^{-2}+p^{-4})
	\\&=\dprod_{p}
\ff{(1-p^{-2})(1-p^{-12})}{(1-p^{-4})(1-p^{-6})}
	\\&=\ff{\z(4)\z(6)}{\z(2)\z(12)}
	\\&=\ff{\p^{4}}{90}\,\ff{\p^{6}}{945}\,\ff{6}{\p^{2}}\,
\ff{638512875}{691\p^{12}}
\vi&=\ff{45045}{691\p^{4}}.
\end{align*}
\end{proof}


\begin{proof}[Proof of \tn{(2)}]
Since 
\[
1-x^{2}+x^{4}-x^{6}=\ff{(1-x^{2})(1-x^{8})}{1-x^{4}},\]
we have 
\begin{align*}
	\dsum_{n:\text{4-free}}\ff{\mu_{3}(n)}{n^{2}}&=
	1+\sum_{
\substack
{0<m_{j}<4
}
}
\ff{(-1)^{m_{1}+\cd+m_{k}}}{(p_{1}^{m_{1}}\cd p_{k}^{m_{k}})^{2}}
	\\&=\dprod_{p}
(1-p^{-2}+p^{-4}-p^{-6})
	\\&=\dprod_{p}
\ff{(1-p^{-2})(1-p^{-8})}{1-p^{-4}}
	\\&=\ff{\z(4)}{\z(2)\z(8)}
	\\&=\ff{\p^{4}}{90}\,\ff{6}{\p^{2}}\,\ff{9450}{\p^{8}}
\vi&=\ff{630}{\p^{6}}.
\end{align*}
\end{proof}


\begin{proof}[Proof of \tn{(3)}]
The idea is quite similar.
From the cyclotomic polynomial
\[
\Phi_{20}(x)=
\ff{(1-x^{2})(1-x^{20})}{(1-x^{4})(1-x^{10})}
=1-x^{2}+x^{4}-x^{6}+x^{8}, 
\]
we obtain 
\begin{align*}
	\di\sum_{
	\substack
	{n:
	\textnormal{5-free}
	}
	}
	\ff{\mu_{5}(n)}{n^{2}}&=
	1+\sum_{
\substack
{0<m_{j}<5
}
}
\ff{(-1)^{m_{1}+\cd+m_{k}}}{(p_{1}^{m_{1}}\cd p_{k}^{m_{k}})^{2}}
	\\&=\prod_{p}(1-p^{-2}+p^{-4}-p^{-6}+p^{-8})
	\\&=\prod_{p}\ff{(1-p^{-2})(1-p^{-20})}{(1-p^{-4})(1-p^{-10})}
	\\&=\ff{\z(4)\z(10)}{\z(2)\z(20)}
	\\&=\ff{\p^{4}}{90}\,
\ff{\p^{10}}{93555}\,
	\ff{6}{\p^{2}}\,\ff{1531329465290625}{174611\p^{20}}
\vi&=\ff{1091215125}{174611\p^{8}}.
\end{align*}
\end{proof}

\subsection{M\"obius inversion of rank 3}

Recall that $\mu_{2}^{-1}=1$. Thus, the equality 
\[
\k{\dsum_{n=1}^{\mug}\ff{\mu(n)}{n^{s}}}^{-1}
=\dsum_{n=1}^{\mug}\ff{1}{n^{s}}
\]
can be regarded as ``M\"obius inversion of rank 2".
Here we consider the case of rank 3.

For the prime factorization $n=p_{1}^{m_{1}}\cd p_{k}^{m_{k}}$ of $n$ into distinct primes, let
\[
m_{*}(n)=|\{j\,|\,m_{j}\equiv 3, \tn{4 mod 6}\}|.
\]

\begin{thm}[M\"obius inversion of rank 3]
\[
\k{\dsum_{n=1}^{\mug}\ff{\mu_{3}(n)}{n^{s}}}
^{-1}=
\sum_{
\substack
{n=p_{1}^{m_{1}}\cd p_{k}^{m_{k}}\\
m_{j}\not\equiv 2, \tn{5 mod }6
}
}
\ff{(-1)^{m_*(n)}}{n^{s}}.
\]\end{thm}

\begin{proof}
We know that 
\[
\dsum_{n=1}^{\mug}\ff{\mu_{3}(n)}{n^{s}}
=
\dsum_{n:\textnormal{3-free}}\ff{\mu_{3}(n)}{n^{s}}=
\prod_{p}(1-(p^{-1})^{s}+(p^{-2})^{s}).
\]
Now the idea is to find the inverse (formal  power series) of $1-x+x^{2}$:
\begin{align*}
	(1-x+x^{2})^{-1}&=\Phi_{6}(x)^{-1}
	\\&=\ff{(1-x^{2})(1-x^{3})}{(1-x)(1-x^{6})}
	\\&=\ff{1+x-x^{3}-x^{4}}{1-x^{6}}
	\\&=\dsum_{i=0}^{\mug}x^{6i}(1+x-x^{3}-x^{4}).
\end{align*}
That is, 
\[
\dsum_{n=1}^{\mug}\ff{\mu_{3}^{-1}(n)}{n^{s}}
=
\dprod_{p}\k{
\dsum_{i=0}^{\mug}p^{-6si}(1+p^{-s}-p^{-3s}-p^{-4s})}.
\]
It follows that
\[
\dsum_{n=1}^{\mug}\ff{\mu_{3}^{-1}(n)}{n^{s}}
=
\dprod_{p}\k{\dsum_{i=0}^{\mug}p^{-6si}(1+p^{-s}-p^{-3s}-p^{-4s})}
=
\sum_{
\substack
{n=p_{1}^{m_{1}}\cd p_{k}^{m_{k}}\\
m_{j}\not\equiv 2, \tn{5 mod }6
}
}
\ff{(-1)^{m_*(n)}}{n^{s}}.
\]

{\renewcommand{\arraystretch}{1.5}
\begin{table}
\caption{$\mu_{3}(n)$ and $\mu_{3}^{-1}(n)$}
\begin{center}
\begin{tabular}{c|ccccccccccccccccccc}
$n$&1&2&3&4&5&6&7&8&9&10&$\cd$\\\hline
$\mu_{3}(n)$&1&$-1$&$-1$&1&$-1$&1&$-1$&0&1&1&$\cd$\\
$\mu_{3}^{-1}(n)$&1&1&1&0&1&1&1&$-1$&0&1&$\cd$\\
\end{tabular}
\end{center}
\end{table}%
}
\end{proof}

\begin{cor}
\[
\mu_{3}^{-1}(p^{m})=
\begin{cases}
	1&m\equiv 0, \tn{1  \,\,mod 6},\\
	0&m\equiv 2, \tn{5  \,\,mod 6},\\
	-1&m\equiv 3, \tn{4  \,\,mod 6}.
\end{cases}
\]
\end{cor}

\begin{cor}
\[
\dsum_{n=1}^{\mug}\ff{\mu_{3}^{-1}(n)}{n^{2}}=\ff{691\p^{4}}{45045}.
\]
\end{cor}
\begin{proof}
This is the inverse of the sum
\[
\k{\dsum_{n=1}^{\mug}\ff{\mu_{3}(n)}{n^{2}}
=}
\dsum_{n:\tn{3-free}}\ff{\mu_{3}(n)}{n^{2}}=
\ff{45045}{691\p^{4}}
\]
proved in Theorem \ref{mth1}.
\end{proof}
Let us put it this way; this gives a new kind of an infinite series expression of $\p$ in terms of $\mu_{3}^{-1}$:
\begin{align*}
	\p&=
	\k{\ff{45045}{691}
	\k{1+\ff{1}{2^{2}}+\ff{1}{3^{2}}+\ff{1}{5^{2}}+\ff{1}{6^{2}}
+\ff{1}{7^{2}}-\ff{1}{8^{2}}+\ff{1}{10^{2}}+\cd}}^{1/4}.
\end{align*}

%
%
%

\subsection{Theorem ($s=3$)}

We evaluated several Dirichlet series at $s=2$ so that many zeta values at even integers $\{\z(2n)\}$ appeared. 
On the one hand, exact values of $\z(2n)$ are known in terms of \eh{Bernoulli numbers}; on the other hand, not much is known on $\{\z(2n+1)\}$.
\begin{itemize}
	\item Ap\'{e}ry \cite{apery} proved that $\z(3)$ is irrational in 1979.
	\item More recently, Zudilin \cite{zudilin} proved that at least  
	one of $\z(5), \z(7), \z(9), \z(11)$ is irrational.
	\item Exact value of any $\z(2n+1)$ is not known.
\end{itemize}

However, our method is helpful for understanding some \eh{relation} of  particular series involving $\{\zeta(2n+1)\}$.
Let us see an example on $s=3$, $\zeta(3)$ and $\zeta(9)$ here.

\begin{lem}\label{lemma}
\hf\begin{quote}
	\begin{enumerate}
		\item $\dsum_{n:\tn{3-free}}\ff{1}{n^{3}}=
		\ff{\z(3)}{\z(9)}$.
		\item 
		$\dsum_{n:\tn{3-free}}\ff{\mu_{3}(n)}{n^{3}}=
		\ff{\z(6)\z(9)}{\z(3)\z(18)}$.
	\end{enumerate}
\end{quote}
\end{lem}

\begin{proof}
Take the cyclotomic polynomials
\begin{align*}
	\Phi_{9}(x)&=\ff{1-x^{9}}{1-x^{3}}=1+x^{3}+x^{6} \tn{\mb{ } and }
	\\\Phi_{18}(x)&=\ff{(1-x^{3})(1-x^{18})}{(1-x^{6})(1-x^{9})}
	=1-x^{3}+x^{6}.
\end{align*}

Then 
\[
\dsum_{n:\tn{3-free}}\ff{1}{n^{3}}=
\dprod_{p}(1+(p^{-1})^{3}+(p^{-2})^{3})
=\dprod_{p}\ff{1-p^{-9}}{1-p^{-3}}=\ff{\z(3)}{\z(9)} \tn{\q and}
\]
\[
\dsum_{n:\tn{3-free}}\ff{\mu_{3}(n)}{n^{3}}=
\dprod_{p}(1-(p^{-1})^{3}+(p^{-2})^{3})
=\dprod_{p}\ff{(1-p^{-3})(1-p^{-18})}{(1-p^{-3})(1-p^{-6})}
=\ff{\z(6)\z(9)}{\z(3)\z(18)}.
\]
\end{proof}

\begin{thm}
$
		\k{\di\sum_{
	\substack
	{n:
	\textnormal{3-free}
	}
	}
	\ff{1}{n^{3}}
}
	\k{\di\sum_{
	\substack
	{n:
	\textnormal{3-free}
	}
	}
	\ff{\mu_{3}(n)}{n^{3}}}
	=\ff{41247931725}{43867\p^{12}}.
	$
\end{thm}

\begin{proof}
Thanks to Lemma \ref{lemma}, the left hand side is 
\[
\ff{\z(6)}{\z(18)}=
\ff{\p^{6}}{945}\,\ff{38979295480125}{43867\p^{18}}
=\ff{41247931725}{43867\p^{12}}
.\]
\end{proof}


\section{Final remarks}

\subsection{Lambert series}

Here, we record some of our results in a little different form. 

For a sequence $a_{n}$ of integers, its \eh{Lambert series} is the formal power series
\[
\dsum_{n=1}^{\mug}\ff{a_{n}x^{n}}{1-x^{n}}.\]
Assume that $a_{n}=f(n)$ for some arithmetic function $f$. It turns out that the coefficient of $x^{N}$ in $\dsum_{n=1}^{\mug}\ff{f(n)x^{n}}{1-x^{n}}$ is 
exactly $\dsum_{d|N}f(d)$, that is, $(f*1)(N)$.
\begin{cor}
For $r\ge 3$ odd or $r=\mug$, we have 
\[
\dsum_{n=1}^{\mug}\ff{\murn x^{n}}{1-x^{n}}=
\sum_{n=1}^{\mug}x^{n^{2}}.
\]
\end{cor}
In particular, this includes
\[
\dsum_{n=1}^{\mug}\ff{\lam(n)x^{n}}{1-x^{n}}=
\sum_{n=1}^{\mug}x^{n^{2}}
=x+x^{4}+x^{9}+x^{16}+\cd
\]
as a special case.

\begin{cor}
For $r$ even, we have 
\[
\dsum_{n=1}^{\mug}\ff{\mur(n) x^{n}}{1-x^{n}}=
\sum_{
\substack
{n:n= p_{1}^{m_{1}}\cd p_{k}^{m_{k}}\\
m_{j}<r, \,\,m_{j}\,\, \tn{ even}
}
}
x^{n}.
\]
\end{cor}

For example, $r=4$, 
\[
\ff{\mu_{4}(1)x}{1-x}+
\ff{\mu_{4}(2)x^{2}}{1-x^{2}}+
\ff{\mu_{4}(3)x^{3}}{1-x^{3}}+
\ff{\mu_{4}(4)x^{4}}{1-x^{4}}+
\ff{\mu_{4}(5)x^{5}}{1-x^{5}}+
\cd
\]
\[=
x+x^{4}+x^{9}+x^{25}+x^{36}+x^{49}
+x^{100}+x^{121}+x^{169}+\cd.\]

\subsection{Future research}

We leave several ideas here for our future research.

\begin{enumerate}
	\item We expect that there are many more results on Dirichlet series  
	\[
\dsum_{n:\textnormal{$r$-free}}\ff{f(n)}{n^{s}}
, \q r, s\ge 2, \q f\in\{\mu_{r}, \mu_{r}*1, \mu_{r}^{-1}, 1\}.
\]
\item 
Suppose a multiplicative arithmetic function $f$ satisfies
\[
f(p^{m})=f(q^{m}) \textnormal{\q for all primes $p, q$.}
\]
The \eh{Bell series} for such $f$ is the formal power series 
\[
B_{f}(x)=\sum_{m=0}^{\mug}f(p^{m})x^{m} 
\]
as $B_{\mu}(x)=1-x$, for instance.
Say $f$ is \eh{cyclotomic} if 
$B_{f}(x)=\Phi_{n}(x)$ for some $n$;
it is \eh{inverse cyclotomic} if 
$B_{f}(x)=\Phi_{n}(x)^{-1}$ for some $n$.
Study a series $\di\sum_{n}\ff{f(n)}{n^{s}}$ for functions of this class.
\item Describe details of $\mu_{r}^{-1}$ for $r\ge 4$.
\end{enumerate}


\np

\end{document}